\documentclass{ifacconf}

\usepackage{natbib}
\usepackage{cuted}
\usepackage{graphicx}
\usepackage{url}
\usepackage{mathtools}
\usepackage{tikz}
\usetikzlibrary{shapes,arrows}
\usetikzlibrary{automata}
\usepackage{amsmath}
\usepackage{verbatim}
\usepackage{listings}
\usepackage{amssymb}
\usepackage{mathrsfs}
\usepackage{enumerate}
\allowdisplaybreaks

\usepackage{tikz}
\usepackage{textcomp}

\tikzstyle{block} = [draw, rectangle,
    minimum height=3em, minimum width=4em,align=center]
\tikzstyle{hblock} = [draw, rectangle,
    minimum height=15em, minimum width=3em,align=center]
\tikzstyle{outerblock} = [draw, rectangle, dashed, inner sep = .85em]
\tikzstyle{outerellipse} = [draw, rectangle, dashed, inner sep = 0.3em]
\tikzstyle{circ} = [draw, circle, 
   node distance=3em ,align=center, fill=gray!0!white]
\tikzstyle{sum} = [draw, circle, radius = .5cm, node distance=1cm]
\tikzstyle{int} = [coordinate]
\tikzstyle{input} = [coordinate]
\tikzstyle{output} = [coordinate]
\tikzstyle{pinstyle} = [pin edge={to-,thin,black}]
\usepackage[most]{tcolorbox}

\definecolor{Ablue}{rgb}{0.16, 0.32, 0.75}
\definecolor{Ayellow}{rgb}{1.0, 0.92, 0.8}
\definecolor{pred}{rgb}{0.8, 0, 0}
\definecolor{pgreen}{rgb}{0, 0.8, 0}
\definecolor{pblue}{rgb}{0, 0, 0.8}

\definecolor{amber}{rgb}{1.0, 0.75, 0.0}

{\theorembodyfont{\slshape}\newtheorem{theorem}{Theorem}}
{\theorembodyfont{\upshape}\newtheorem{proposition}[theorem]{Proposition}}  %
{\theorembodyfont{\slshape}\newtheorem{lemma}[theorem]{Lemma}}
{\theorembodyfont{\slshape}}           %
{\theorembodyfont{\slshape}}
{\theorembodyfont{\upshape}\newtheorem{corollary}[theorem]{Corollary}}  %
{\theorembodyfont{\upshape}\newtheorem{definition}[theorem]{Definition}}
{\theorembodyfont{\upshape}}
{\theorembodyfont{\upshape}\newtheorem{remark}[theorem]{Remark}}
{\theorembodyfont{\upshape}}
{\theorembodyfont{\upshape}\newtheorem{example}[theorem]{Example}}
{\theorembodyfont{\upshape}}

\usepackage{pgfplots}
\pgfplotsset{compat=1.11}
\usepgfplotslibrary{fillbetween}
\pgfdeclarelayer{bg}    
\pgfsetlayers{bg,main}  

\begin{document}
\begin{frontmatter}
\title{Informativity conditions for data-driven control based on input-state data and polyhedral cross-covariance noise bounds}\thanks{This work has received funding from the European Research Council (ERC), Advanced Research Grant SYSDYNET, under the European Union's Horizon 2020 research and innovation programme (Grant Agreement No. 694504).}

\author{Tom R.V. Steentjes, Mircea Lazar, Paul M.J. Van den Hof}
\address{Department of Electrical Engineering, Eindhoven University of Technology, The Netherlands (e-mails:  t.r.v.steentjes@tue.nl, m.lazar@tue.nl, p.m.j.vandenhof@tue.nl)}

\begin{abstract}
Modeling and control of dynamical systems rely on measured data, which contains information about the system. Finite data measurements typically lead to a set of system models that are unfalsified, i.e., that explain the data. The problem of data-informativity for stabilization or control with quadratic performance is concerned with the existence of a controller that stabilizes all unfalsified systems or achieves a desired quadratic performance. Recent results in the literature provide informativity conditions for control based on input-state data and ellipsoidal noise bounds, such as energy or magnitude bounds. In this paper, we consider informativity of input-state data for control where noise bounds are defined through the cross-covariance of the noise with respect to an instrumental variable; bounds that were introduced originally as a noise characterization in parameter bounding identification. The considered cross-covariance bounds are defined by a finite number of hyperplanes, which induce a (possibly unbounded) polyhedral set of unfalsified systems. We provide informativity conditions for input-state data with polyhedral cross-covariance bounds for stabilization and $\mathcal{H}_2$/$\mathcal{H}_\infty$ control through vertex/half-space representations of the polyhedral set of unfalsified systems.
\end{abstract}

\end{frontmatter}
\section{Introduction}
Models of dynamical systems play a key role in the synthesis of controllers. Typically, these models are not available, however, and have to be derived from data and prior knowledge from first-principles modelling. Estimating dynamical models from measurement data is considered in the field of system identification~\citep{ljung1999}. Based on models identified from data, controllers can be indirectly synthesized via model-based control methods through the certainty equivalence principle~\citep{hou2013}. This is also referred to as indirect data-driven control in the literature. Taking the control objective into account in the identification can lead to models that are especially fit control design; a topic that has been extensively studied in the field of identification for control~\citep{vandenhofetal1995}.

The modelling step in data-driven control may be circumvented to synthesize a controller directly based on the data. Methods for direct data-driven control include adaptive control methods, virtual reference feedback tuning~\citep{campi2002}, iterative feedback tuning~
\citep{hjalmarsson98}, and optimal controller identification~\citep{campestrini2017}, see e.g.,~\citep{hou2013} for an overview of methods for data-driven control. A commonality of the mentioned data-driven control methods is that persistently exciting data are required, i.e., the data are in fact informative enough for system identification.

Even if data are not informative enough for system identification, data can still be informative enough for controller design. The pioneering work~\citep{waarde2020tac} introduced a framework for analyzing informativity of data for system-theoretic properties and controller design. In particular, necessary and sufficient conditions for informativity of noiseless data for controller design were developed in~\citep{waarde2020tac}, which can hold even if the data are not informative for system identification.

In general, process noise will be present, but prior knowledge on the noise, if available, can be taken into account in the informativity analysis. This problem has recently received considerable attention in the literature~\citep{berberich2020}, \citep{vanwaarde2020noisytac}, \citep{bisoffi2021}, \citep{vanwaarde2021finsler}, \citep{steentjes2021cdc}. The prior knowledge considered in the aforementioned literature is typically represented by a quadratic bound on the noise sequence, which includes bounds on the energy and magnitude of the noise process. Prior knowledge of the noise in the form of ellipsoidal bounds on the sample cross-covariance have been considered in~\citep{steentjes2022qccb}. Sample cross-covariance bounds were introduced in~\citep{hakvoort95} as an alternative to magnitude bounds in parameter bounding identification, given its overly conservative noise characterization.

In this paper, we consider informativity of input-state data for controller design in the presence of noise satisfying \emph{polyhedral} cross-covariance bounds. This prior knowledge combined with measurement data leads to sets of feasible system matrices that are intersections of halfspaces and therefore (possibly unbounded) polyhedra. We show how convexity of the sets of feasible system matrices and stability/performance criteria lead to data-based linear matrix inequalities (LMIs) that are necessary and sufficient for quadratic stabilization, $\mathcal{H}_\infty$ and $\mathcal{H}_2$ control in the case the polyhedron is bounded. The technique of using the convexity of polytopes for obtaining a finite set of controller synthesis LMIs is well known in robust control, e.g. for stabilization of systems with polytopic uncertainties, cf.~\citep{kothare1996}, \cite[Chapter~5]{schererweilandLMI}. When the set of feasible system matrices is unbounded, there is no correspondent from robust control for systems with polytopic uncertainty. An unbounded set of feasible systems implies that data are not informative for system identification in the case of noise-free data, cf. \citep[Example~19]{waarde2020tac}, and is therefore particularly interesting for informativity analysis.  We provide preliminary results for data informativity for stabilization, in the case of noisy data that lead to a unbounded set of feasible systems.

\section{Polyhedral cross-covariance bounds}
In this paper, we consider the data-informativity for a class of linear systems that is affected by a noise signal $e(t)$:
\begin{align} \label{eq:sysclass}
x(t+1)=Ax(t)+Bu(t)+e(t),
\end{align}
with state dimension $n$ and input dimension $m$.

The true system is represented by the pair $(A_0,B_0)$. State and input data generated by the true system are collected in the matrices
\begin{align*}
X := [x(0)\ \cdots\ x(N)],\quad U_-:=[u(0)\  \cdots\ u(N-1)].
\end{align*}
By defining
\begin{align*}
X_+ &:= [x(1)\ \cdots\ x(N)],\quad X_-:=[x(0)\ \cdots\ x(N-1)],\\
 E_-&:=[e(0)\ \cdots\ e(N-1)],
\end{align*}
we clearly have
\begin{align} \label{eq:dateq0}
X_+=A_0X_-+B_0U_-+E_-.
\end{align}
In case the noise is \emph{measured}, the set of systems that is consistent with the data $(U_-,X)$ is
\begin{align*}
\Sigma_{(U_-,X,E_-)} =\{(A,B)\,|\, X_+=AX_-+BU_-+E_-\}.
\end{align*}
When the data are informative for system identification, as defined in \citep{waarde2020tac}, the set of feasible system is a singleton $\Sigma_{(U_-,X,E_-)}=\{(A_0,B_0)\}$. This is equivalent with $\operatorname{col} (X_-, U_-)$ having full rank. In the case the data are not informative, the set $\Sigma_{(U_-,X,E_-)}$ is not a singleton, but becomes a line or hyperplane. Even if the data are not informative for system identification, the data can still be informative for other properties, such as stability or feedback stabilization, cf. \citep{waarde2020tac}.

Let $e=:\operatorname{col}(e_1,\dots,e_n)$ and consider that each noise channel $e_j$, $j=1,\dots,n$, is not measured, i.e., $E_j^-$ is unknown, but that $e_j(t)$ satisfies the bounds
\begin{align} \label{eq:covbndmult}
c_{ij}^l\leq \frac{1}{\sqrt{N}}\sum_{t=0}^{N-1} r_i(t)e_j(t)\leq c_{ij}^u,\quad i=1,\dots,M,
\end{align}
where $r_i$ are signals that are chosen, typically as a (delayed version of) state or input signal, and $c_{ij}^l$, $c_{ij}^u$ are specified bounds. Notice that we specify $M$ upper and lower bounds for each noise channel $j\in\{1,\dots,n\}$, and that the instrumental variables $r_i$, $i\in\{1,\dots,M\}$, are common for all noise channels $j\in\{1,\dots,n\}$. The bounds in \eqref{eq:covbndmult} are satisfied for all $j$ if and only if 
\begin{align*}
E_-=\operatorname{col}(E_1^-,\dots,E_n^-)\in \mathcal{E}_R,
\end{align*}
where
\begin{align*}
\mathcal{E}_R&:=\{E\,|\, C_l\leq \frac{1}{\sqrt{N}}E_-R_-^\top\leq C_u\}\\
&=\{E\,|\, C_l\leq \frac{1}{\sqrt{N}}\sum_{t=0}^{N-1} e(t)r(t)^\top\leq C_u\}.
\end{align*}
with $R_-:=\operatorname{col} (R_1^-,\dots R_M^-)$, $R_i^-:=[r_i(0)\cdots r_i(N-1)]$, and with $c_{ij}^l$ and $c_{ij}^u$ the $(i,j)$-th entry of $C_l$ and $C_u$, respectively. The inequalities defining $\mathcal{E}_R$ are thus \emph{entry-wise} inequalities.

\begin{remark}
Noise bounds of the type~\eqref{eq:covbndmult} define upper and lower bounds on the sample cross-covariance of the noise $e$ and an instrumental variable $r$. These bounds were introduced in~\citep{hakvoort95} for parameter bounding identification. An `ellipsoidal' version of these bounds, i.e., a bound on $E_-R_-^\top R_-E_-^\top$ in the terms of the partial order on positive semi-definite matrices, has been considered in~\citep{steentjes2022qccb} for analyzing informativity for control. The difference in prior knowledge on the noise has two implications: (i) the bounds~\eqref{eq:covbndmult} allow a component-wise specification of bounds on the cross-covariance compared to ellipsoidal bounds, and (ii) incorporating this ``polyhedral'' (possibly unbounded) prior knowledge on the noise in the informativity analysis requires a fundamentally different approach compared with the application of the matrix S-lemma~\citep{vanwaarde2020noisytac} used in~\citep{steentjes2022qccb}, as will be discussed in Section~\ref{sec:inf}.
\end{remark}
\begin{remark}
Guidelines in the literature recommend choosing an instrumental variable $r$ that is correlated with the input $u$, but uncorrelated with the noise $e$~\citep{hakvoort95}. We refer to~\citep{hakvoort95} for more information on choosing $r$ and estimating the bounds~\eqref{eq:covbndmult} from data.
\end{remark}
The bounds on the cross-covariance between the noise channels and the instrumental signals induce a restriction on the pairs $(A,B)$ that satisfy the data equation
\begin{align} \label{eq:dateq}
X_+=AX_-+BU_-+E_-.
\end{align}
All systems that explain the data $(U_-,X)$ for some $E_-\in\mathcal{E}_R$ are collected in the set $\Sigma_{(U_-,X)}^R$:
\begin{align*}
\Sigma_{(U_-,X)}^R:=\{(A,B)\,|\, \exists E_-\in\mathcal{E}_R \text{ such that } \eqref{eq:dateq} \text{ holds}\}.
\end{align*}

The following proposition provides a parametrization for the set of feasible systems with cross-covariance bounds.
\begin{proposition} $\Sigma_{(U_-,X)}^R=\{(A,B)\,|\, \eqref{eq:lccbparam} \text{ holds}\}$,
where
\begin{align} \label{eq:lccbparam}
\sqrt{N}C_l\leq X_+R_-^\top-\begin{bmatrix} A & B \end{bmatrix} \begin{bmatrix}X_-R_-^\top\\ U_-R_-^\top \end{bmatrix}\leq \sqrt{N}C_u
\end{align}
\end{proposition}
\begin{pf}
The set of feasible system matrices is
\begin{align} \label{eq:sigmarhalfspace}
\Sigma_{(U_-,X)}^R&=\{(A,B)\,|\, C_l\leq \frac{1}{\sqrt{N}}\sum_{t=0}^{N-1}e(t)r(t)^\top\leq C_u\}\\
&=\{(A,B)\,|\, C^l\leq R_{er}^{N-}\leq C_u\},
\end{align}
where
\begin{align*}
R_{er}^{N-}&:=\frac{1}{\sqrt{N}}\sum_{t=0}^{N-1}e(t)r(t)^\top\\
&= \frac{1}{\sqrt{N}}\sum_{t=0}^{N-1}\left(x(t+1)-Ax(t)-Bu(t)\right)r(t)^\top\\
&=\frac{1}{\sqrt{N}}X_+R_-^\top-A\frac{1}{\sqrt{N}}X_-R_-^\top-B\frac{1}{\sqrt{N}}U_-R_-^\top.
\end{align*}
Hence, the feasible set of systems is 
\begin{align*}
\Sigma_{(U_-,X)}^R=\{(A,B)\,|\, \eqref{eq:lccbparam} \text{ holds}\},
\end{align*}
which completes the proof.\hfill $\square$
\end{pf}

It can be shown that $\Sigma_{(U_-,X)}^R$ is an intersection of half spaces, by observing that
\begin{align*}
\Sigma_{(U_-,X)}^R=\Sigma_{(U_-,X)}^{R_1} \cap\ \cdots\ \cap\Sigma_{(U_-,X)}^{R_M}=\bigcap_{i=1}^M \Sigma_{(U_-,X)}^{R_i},
\end{align*}
where, denoting the $i$-th column of $C_l$ ($C_u$) by $c_i^l$ ($c_i^u$),
\begin{align*}
\Sigma_{(U_-,X)}^{R_i}=\{(A,B)\,|\,c_i^l\leq R_{xr_i}^{N+}-\begin{bmatrix}
A & B
\end{bmatrix}\begin{bmatrix}
R_{xr_i}^{N-}\\ R_{ur_i}^{N-}
\end{bmatrix}\leq c_i^u\},
\end{align*}
with $R_{xr_i}^{N+}=\frac{1}{\sqrt{N}}X_+(R_i^-)^\top$, $R_{xr_i}^{N-}=\frac{1}{\sqrt{N}}X_-(R_i^-)^\top$ and $R_{ur_i}^{N+}=\frac{1}{\sqrt{N}}U_-(R_i^-)^\top$.
Hence, the set of feasible subsystems is either an intersection of halfspaces and unbounded (called an $\mathcal{H}$-polyhedron) or it is a bounded polyhedron (called $\mathcal{V}$-polytope). Another way to see that $\Sigma_{(U_-,X)}^R$ is an intersection of halfspaces, is to vectorize the inequalities:
\begin{align*}
\Sigma_{(U_-,X)}^R=\{(A,B)\,|\, &\operatorname{vec}(C_l)\leq \operatorname{vec}(R_{xr}^{N+})-\!\left(\!\begin{bmatrix}R_{xr}^{N-}\\ R_{ur}^{N-} \end{bmatrix}^\top\!\!\!\!\otimes I_n\!\right)\\
&\times \operatorname{vec}\left(\begin{bmatrix} A & B \end{bmatrix}\right) \leq \operatorname{vec}(C_u)\}.
\end{align*}
\begin{lemma} \label{lem:kerxr}
The set of feasible systems $\Sigma_{(U_-,X)}^R$ is bounded if and only if
\begin{align} \label{eq:kerxr}
\operatorname{ker} \begin{bmatrix}
X_-R_-^\top\\ U_-R_-^\top
\end{bmatrix}^
\top=\{0\}.
\end{align}
\end{lemma}
\begin{pf}
First, we note that $\Sigma_{(U_-,X)}^R$ is not empty. A non-empty polyhedron
\begin{align*}
\Sigma_{(U_-,X)}^R=\{(A,B) \,|\, M\operatorname{vec}(\begin{bmatrix}
A & B
\end{bmatrix})\leq c\}
\end{align*}
is unbounded if and only if there exists $v\neq 0$ so that $Mv\leq 0$. With
\begin{align*}
M:=\begin{bmatrix}
-\left(\begin{bmatrix}R_{xr}^{N-}\\ R_{ur}^{N-} \end{bmatrix}^\top\otimes I_n\right)\\
\left(\begin{bmatrix}R_{xr}^{N-}\\ R_{ur}^{N-} \end{bmatrix}^\top\otimes I_n\right)
\end{bmatrix},
\end{align*}
we observe that $Mv\leq 0$ if and only if $Mv=0$. Hence, $\Sigma_{(U_-,X)}^R$ is unbounded if and only if
\begin{align*}
\operatorname{ker}\left(\begin{bmatrix}R_{xr}^{N-}\\ R_{ur}^{N-} \end{bmatrix}^\top\otimes I_n\right)\neq \{0\} \quad \Leftrightarrow\quad \operatorname{ker}\begin{bmatrix}R_{xr}^{N-}\\ R_{ur}^{N-} \end{bmatrix}^\top\neq \{0\}.
\end{align*}
We conclude that $\Sigma_{(U_-,X)}^R$ is bounded if and only if $\operatorname{ker}\begin{bmatrix}R_{xr}^{N-}\\ R_{ur}^{N-} \end{bmatrix}^\top= \{0\}$, which concludes the proof.\hfill $\square$
\end{pf}

\begin{figure}[!t]
\centering
\begin{tikzpicture}
  \begin{axis}[samples = 50, grid=both,ymin=-5,ymax=5,xmax=5,xmin=-5,xticklabel=\empty,yticklabel=\empty,
               minor tick num=1,axis x line = bottom, axis y line = left, xlabel=$A$,ylabel=$B$,label style = {at={(ticklabel cs:0.5)}}]
               
%
               \addplot[green, thick, name path = U] (x,-x+3);
               \addplot[green, thick, name path = L] (x,-x+1);
               \addplot[green!30] fill between[of=L and U];

               \addplot[purple, thick, name path = Ue] (x,.5*x+1);
               \addplot[purple, thick, name path = Le] (x,.5*x-1);
               \addplot[purple!30] fill between[of=Le and Ue];
               
              \addplot[orange!40] fill between[of=L and Ue, soft clip = {domain=0:2/1.5}];
              \addplot[orange!40] fill between[of=Le and U, soft clip = {domain=2/1.5:4/1.5}];
                                                                                  
  \end{axis}
\end{tikzpicture}
\caption{Illustration of the set $\Sigma_{(U_-,X)}^R$ for $M=1$ (green) and for $M>1$ (orange).}
\label{fig:ccbM>1}
\end{figure}
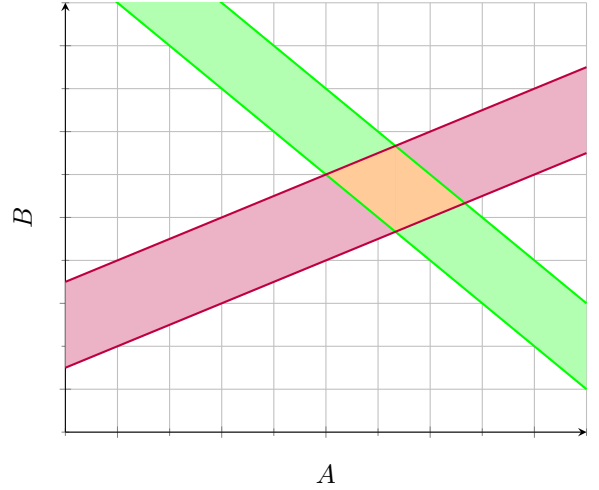

\begin{remark}
The condition for boundedness of $\Sigma_{(U_-,X)}^R$ is equivalent with the matrix $\operatorname{row} (R_-X_-^\top, R_- U_-^\top)$ having full column rank. A necessary condition for the rank of this matrix being full, is to have enough instrumental signals. More precisely, a necessary condition for boundedness is that $M\geq n+m$, where we recall that $n$ and $m$ are the state and input dimension, respectively, and $M$ is the dimension of the instrumental signal $r$. For the scalar case $n=m=1$, an unbounded set $\Sigma_{(U_-,X)}^R$ is obtained for $M=1$, as illustrated in Figure~\ref{fig:ccbM>1} in green. With $M>1$ the rank condition can be satisfied (no redundant inequalities) and a polytope is obtained, as illustrated in Figure~\ref{fig:ccbM>1} in orange.
\end{remark}

\section{Informativity for feedback stabilization} \label{sec:inf}
Consider the problem of stabilizing the `true' system $(A_0,B_0)$ using the data $(U_-,X)$. We define the set of systems that are stabilized\footnote{A matrix is called stable if all its eigenvalues are in the open unit disk.} by $K$ as
\begin{align*}
\Sigma_K:=\{(A,B)\,|\, A+BK\text{ is stable}\}.
\end{align*}
In line with \cite[Definition~12]{waarde2020tac}, we consider the following definition for informativity for stabilization by state feedback.
\begin{definition} \label{def:infstab}
The data $(U_-,X)$ are said to be informative for stabilization by state feedback if there exists a feedback gain $K$ so that 
\begin{align*}
\Sigma_{(U_-,X)}^R\subseteq \Sigma_K.
\end{align*}
\end{definition}
In other words, if there exists a $K$ so that for every system $(A,B)$ in $\Sigma_{(U_-,X)}^R$, $A+BK$ is stable, then the data are informative for feedback stabilization.
\begin{definition} \label{def:infqstab}
The data $(U_-,X)$ are said to be informative for quadratic stabilization by state feedback if there exist a $K$ and $P\succ 0$ so that 
\begin{align} \label{eq:infqstab}
\Sigma_{(U_-,X)}^R\subseteq\! \{(A,B)\,|\, (A+BK) P (A+BK)^\top\!\!-\!P\!\prec 0\}.\!
\end{align}
\end{definition}
Notice the difference: the data are informative for quadratic stabilization if there exists a common pair $(K,P)$, with $P\succ 0$, such that the inclusion in Definition~\ref{def:infqstab} holds, while the data are informative for stabilization if there is a common $K$ so that $\Sigma_{(U_-,X)}^R\subseteq \Sigma_K$. Hence, the data $(U_-,X)$ are informative for stabilization by state feedback if the data $(U_-,X)$ are informative for \emph{quadratic} stabilization by state feedback, but the reverse implication is not true, in general.

\subsection{$\Sigma_{(U_-,X)}^R$ is an unbounded polyhedron} \label{sec:unbnd}
We consider here the scalar case, i.e., $m=n=1$. In the case that there is one instrumental signal $r=r_1$, the set $\Sigma_{(U_-,X)}^R$ is described by two linear inequalities
\begin{align*}
\begin{bmatrix}
A & B
\end{bmatrix}\begin{bmatrix}
R_{xr}^{N-}\\ R_{ur}^{N-}
\end{bmatrix}\leq R_{xr}^{N+}-c^l,\ \begin{bmatrix}
A & B
\end{bmatrix}\begin{bmatrix}
R_{xr}^{N-}\\ R_{ur}^{N-}
\end{bmatrix}\geq R_{xr}^{N+}-c^u.
\end{align*}
We observe that $\Sigma_{(U_-,X)}^R$ is the intersection of two closed half-spaces. The following result states that a sufficient condition for data informativity for stabilization, is the existence of a $K$ that stabilizes all systems on the ``boundaries'', i.e., the defining hyperplanes of $\Sigma_{(U_-,X)}^R$.

\begin{proposition} \label{prop:stabunbounded}
Let $R_{xr}^{N-}$ be non-zero and let there exist $(R_{xr}^{N-})^\dagger$ such that\footnote{Note that in this case ($n=1$), $(R_{xr}^{N-})^\dagger$ is a scalar and is unique.
} $R_{xr}^{N-}(R_{xr}^{N-})^\dagger=1$ and
\begin{align*}
(R_{xr}^{N+}-c^l)(R_{xr}^{N-})^\dagger \quad \text{ and } \quad (R_{xr}^{N+}-c^u)(R_{xr}^{N-})^\dagger
\end{align*}
are stable. Then the data $(U_-,X)$ are informative for stabilization by state feedback. Moreover, $K$ is such that $\Sigma_{(U_-,X)}^R\subseteq \Sigma_K$ if $K=R_{ur}^{N-}(R_{xr}^{N-})^\dagger$, with $(R_{xr}^{N-})^\dagger$ as described above.
\end{proposition}
\begin{pf}
Let $(R_{xr}^{N-})^\dagger$ be non-zero and such that
\begin{align*}
(R_{xr}^{N+}-c^l)(R_{xr}^{N-})^\dagger \quad \text{ and } \quad (R_{xr}^{N+}-c^u)(R_{xr}^{N-})^\dagger
\end{align*}
are stable. We will first show that
\begin{align*}
-1<\begin{bmatrix}
A & B
\end{bmatrix}\begin{bmatrix}
R_{xr}^{N-}\\ R_{ur}^{N-}
\end{bmatrix}(R_{xr}^{N-})^\dagger<1.
\end{align*}
Consider the case that $(R_{xr}^{N-})^\dagger$ is positive. Then
\begin{align*}
\begin{bmatrix}
A & B
\end{bmatrix}\begin{bmatrix}
R_{xr}^{N-}\\ R_{ur}^{N-}
\end{bmatrix}(R_{xr}^{N-})^\dagger&\leq (R_{xr}^{N+}-c^l)(R_{xr}^{N-})^\dagger,\\
\begin{bmatrix}
A & B
\end{bmatrix}\begin{bmatrix}
R_{xr}^{N-}\\ R_{ur}^{N-}
\end{bmatrix}(R_{xr}^{N-})^\dagger&\geq (R_{xr}^{N+}-c^u)(R_{xr}^{N-})^\dagger.
\end{align*}
Furthermore, $R_{xr}^{N+}-c^u\leq R_{xr}^{N+}-c^l$ implies that
\begin{align*}
-1 < (R_{xr}^{N+}-c^u)(R_{xr}^{N-})^\dagger \leq (R_{xr}^{N+}-c^l)(R_{xr}^{N-})^\dagger <1.
\end{align*}
Hence, any $(A,B)\in\Sigma_{(U_-,X)}$ satisfies
\begin{align} \label{eq:ABbnd}
-1<\begin{bmatrix}
A & B
\end{bmatrix}\begin{bmatrix}
R_{xr}^{N-}\\ R_{ur}^{N-}
\end{bmatrix}(R_{xr}^{N-})^\dagger<1.
\end{align}
Similarly, if $(R_{xr}^{N-})^\dagger$ is negative, then $R_{xr}^{N+}-c^u\leq R_{xr}^{N+}-c^l$ implies that
\begin{align*}
-1 < (R_{xr}^{N+}-c^l)(R_{xr}^{N-})^\dagger \leq (R_{xr}^{N+}-c^u)(R_{xr}^{N-})^\dagger <1,
\end{align*}
and we again find that \eqref{eq:ABbnd} for any $(A,B)\in\Sigma_{(U_-,X)}^R$.

Now, define $K:=R_{ur}^{N-}(R_{xr}^{N-})^\dagger$ to observe that $-1 < A+BK <1$ for any $(A,B)\in \Sigma_{(U_-,X)}^R$. Hence, there exists a $K$ so that $\Sigma_{(U_-,X)}^R\subseteq \Sigma_K$, which completes the proof.\hfill $\square$
\end{pf}
\begin{example} \label{example:stabunbnd}
Consider that data $X=[0\ 1.2\ 3\ 4.1\ 4.25]$, $U_-=[1\ 1\ -0.5\ -2]$ have been collected from a system with system matrices $A_0=1.5$ and $B_0=1$. The corresponding noise $E_-=[0.2\ 0.2\ 0.1\ 0.1]$ is unknown, but satisfies $E_-\in\mathcal{E}_R$ for $R_-=U_-$ with $c_u=-c_l=0.25$. For this example, $(R_{xr}^{N+}-c^l)(R_{xr}^{N-})^\dagger=0.6882$ and $(R_{xr}^{N+}-c^u)(R_{xr}^{N-})^\dagger=0.8059$, hence the data are informative for stabilization by state feedback by Proposition~\ref{prop:stabunbounded} and $K=R_{ur}^{N-}(R_{xr}^{N-})^\dagger=-0.7353$ is indeed such that $A_0+B_0K$ is stable.
\end{example}

Alternatively, the sufficient conditions for the data $(U_-,X)$ to be informative for feedback stabilization can be stated in terms of linear matrix inequalities.
\begin{proposition}
Let there exist a $\Theta$ satisfying $R_{xr}^{N-}\Theta = (R_{xr}^{N-} \Theta)^\top$ so that \begin{align}
\begin{bmatrix}
R_{xr}^{N-} \Theta & (R_{xr}^{N+}-c^l)\Theta\\
\Theta^\top(R_{xr}^{N+}-c^l)^\top & R_{xr}^{N-} \Theta
\end{bmatrix}&\succ 0 \quad \text{ and }\label{eq:thetalmil}\\
\begin{bmatrix}
R_{xr}^{N-} \Theta & (R_{xr}^{N+}-c^u)\Theta\\
\Theta^\top(R_{xr}^{N+}-c^u)^\top & R_{xr}^{N-} \Theta
\end{bmatrix}&\succ 0. \label{eq:thetalmiu}
\end{align}
Then the data $(U_-,X)$ are informative for stabilization by state feedback. Moreover, $K$ is such that $\Sigma_{(U_-,X)}^R\subseteq \Sigma_K$ if $K=R_{ur}^{N-}\Theta((R_{xr}^{N-})^\dagger\Theta)^{-1}$.
\end{proposition}
\begin{pf}
The inequalities in~\eqref{eq:thetalmil}-\eqref{eq:thetalmiu} imply that $R_{xr}^{N-}\Theta$ is positive definite and that
\begin{align*}
[(R_{xr}^{N+}-c^l)\Theta(R_{xr}^{N-} \Theta)^{-1}](R_{xr}^{N-} \Theta)[\star]^\top-R_{xr}^{N-}\Theta &<0 \quad \text{and}\\
[(R_{xr}^{N+}-c^u)\Theta(R_{xr}^{N-} \Theta)^{-1}](R_{xr}^{N-} \Theta)[\star]^\top-R_{xr}^{N-}\Theta &<0.
\end{align*}
Hence, $(R_{xr}^{N+}-c^l)\Theta(R_{xr}^{N-} \Theta)^{-1}$ and $(R_{xr}^{N+}-c^u)\Theta(R_{xr}^{N-} \Theta)^{-1}$ are stable. That is, there exists a right inverse $(R_{xr}^{N-})^\dagger:=\Theta (R_{xr}^{N-}\Theta)^{-1}$ so that$(R_{xr}^{N+}-c^l)\Theta(R_{xr}^{N-})^\dagger$ and $(R_{xr}^{N+}-c^u)\Theta(R_{xr}^{N-})^\dagger$ are stable. Therefore, it follows by Proposition~\ref{prop:stabunbounded} that the data $(U_-,X)$ are informative for stabilization by state feedback.\hfill $\square$
\end{pf}

\subsection{$\Sigma_{(U_-,X)}^R$ is a bounded polyhedron} \label{sec:bnd}
By Lemma~\ref{lem:kerxr}, we observe that $\Sigma_{(U_-,X)}^R$ is a convex polytope with a finite number of vertices $\sigma_{(U_-,X)}^i$, $i=1,\dots, L$, if the data $(U_-,X)$ and instrumental signals $R_-$ satisfy \eqref{eq:kerxr}. In the scalar case, for example, the set $\Sigma_{(U_-,X)}^R$ is then described by $L=4$ vertices with $M=2$ instrumental variables, as depicted in Figure~\ref{fig:ccbM>1}.

By Definition~\ref{def:infstab}, the data $(U_-,X)$ are informative for stabilization by state feedback if there exists a $K$ so that $A+BK$ is stable for all $(A,B)\in\Sigma_{(U_-,X)}^R$. If~\eqref{eq:kerxr} holds true, then $\Sigma_{(U_-,X)}^R=\operatorname{conv}\{\sigma_{(U_-,X)}^1,\dots,\sigma_{(U_-,X)}^L\}$. The following lemma allows us to verify stability conditions for al matrices $(A,B)$ that are compatible with the data, by verifying the conditions at the extreme points of $\Sigma_{(U_-,X)}^R$.
\begin{lemma} \label{lem:uniformbound}
Let $\Gamma\in\mathbb{S}^{n\times n}$,\footnote{$\mathbb{S}^{n\times n}$ denotes the set of $n\times n$ symmetric matrices with real entries.} let $\mathcal{S}_0$ be a set and let $F:\mathcal{S}\to \mathbb{S}^{n\times n}$ be a function with domain $\mathcal{S}=\operatorname{conv} \mathcal{S}_0$. Then $F(x)\prec \Gamma$ for all $x\in\mathcal{S}$ if and only if $F(x)\prec \Gamma$ for all $x\in\mathcal{S}_0$.
\end{lemma}
\begin{pf}
The assertion is a strict version of the assertion in~\citep[Proposition~1.14]{schererweilandLMI}. The proof follows \emph{mutatis mutandis} by the proof of \citep[Proposition~1.14]{schererweilandLMI}.\hfill $\square$
\end{pf}
Now, given the (known) vertices $\sigma_{(U_-,X)}^i$, $i=1,\dots, L$, the problem of verifying informativity for stabilization can be reduced to verifying the stability condition at the extreme points of $\Sigma_{(U_-,X)}^R$, as shown by the following result:
\begin{proposition} \label{prop:infstab}
Let \eqref{eq:kerxr} hold. The data $(U_-,X)$ are informative for quadratic stabilization by state feedback if and only if there exist $K$ and $P$ so that $P\succ 0$ and
\begin{align} \label{eq:infstab}
\begin{bmatrix}
I\\ K
\end{bmatrix}^\top (\sigma_{(U_-,X)}^i)^\top P \sigma_{(U_-,X)}^i \begin{bmatrix}
I\\ K
\end{bmatrix}- P\prec 0, \quad i=1,\dots, L.
\end{align}
\end{proposition}
\begin{pf}
Consider the function $F:\Sigma_{(U_-,X)}^R\to \mathbb{S}^{n\times n}$, defined by $F(\sigma):=\operatorname{col}(I,K)^\top\sigma^\top P \sigma \operatorname{col}(I,K)$. Since $\Sigma_{(U_-,X)}^R$ is convex and $P\succ 0$, we infer that $F$ is a convex function. Hence, by Lemma~\ref{lem:uniformbound}, $F(\sigma)\prec P$ for all $\sigma\in \Sigma_{(U_-,X)}^R$ if and only if $F(\sigma)\prec P$ for all $\sigma\in\{\sigma_{(U_-,X)}^1,\dots, \sigma_{(U_-,X)}^L\}$. This proves the assertion.\hfill $\square$
\end{pf}
We note that the conditions in Proposition~\ref{prop:infstab} are not linear with respect to $K$ and $P$. The application of the Schur complement yields conditions equivalent to \eqref{eq:infstab} that are LMIs:
\begin{corollary} \label{cor:quadstab}
Let \eqref{eq:kerxr} hold. The data $(U_-,X)$ are informative for quadratic stabilization by state feedback if and only if there exist $Y$ and $M$ so that
\begin{align} \label{eq:infstablin}
\begin{bmatrix}
Y	&Z^\top (\sigma_{(U_-,X)}^i)^\top\\ \sigma_{(U_-,X)}^i Z &Y
\end{bmatrix}\succ 0, \quad i=1,\dots, L,
\end{align}
with $Z:=\operatorname{col}(Y,M)$. Moreover, $K$ is such that $\Sigma_{(U_-,X)}^R\subseteq \Sigma_K$ if $K=MY^{-1}$.
\end{corollary}
\begin{pf}
By the Schur complement, the existence of $K$ and $P\succ 0$ such that~\eqref{eq:infqstab} is equivalent with
\begin{align*}
\exists K,P\quad \text{ such that }\begin{bmatrix}
P & (A+BK)^\top\\ A+BK & P^{-1}
\end{bmatrix}\succ 0
\end{align*}
for all $(A,B)\in \Sigma_{(U_-,X)}^R$. Define $Y:=P^{-1}$ and $M:=KP^{-1}$ and perform a congruence transformation with $\operatorname{diag}(Y,I)$ to obtain
\begin{align*}
\exists Y,M\quad \text{ such that }\begin{bmatrix}
Y & (AY+BM)^\top \\ AY+BM & Y
\end{bmatrix}\succ 0
\end{align*}
for all$(A,B)\in \Sigma_{(U_-,X)}^R$. By Lemma~\ref{lem:uniformbound}, we find that this is equivalent with~\eqref{eq:infstablin}, which proves the assertion.\hfill $\square$
\end{pf}

\begin{corollary}
Let \eqref{eq:kerxr} hold. The data $(U_-,X)$ are informative for stabilization by state feedback if one (and therefore all) of the following equivalent statements holds:
\begin{itemize}
\item the data $(U_-,X)$ are informative for quadratic stabilization by state feedback,
\item there exist $K$ and $P$ so that $P\succ 0$ and \eqref{eq:infstab} are satisfied,
\item there exist $Y$ and $M$ so that \eqref{eq:infstablin} is satisfied.
\end{itemize}
\end{corollary}

\begin{figure}
\centering
\includegraphics[width=3.5in]{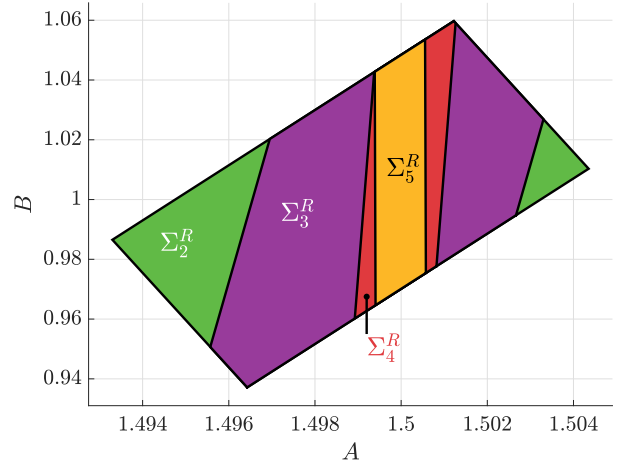}
\caption{Feasible sets of systems $\Sigma_{(U_-,X)}^R$ obtained in Example~\ref{example:stabbnd} with different choices of $R_-$ for $M\in\{2,3,4,5\}$.}
\label{fig:lccbmtns}
\end{figure}

\begin{example} \label{example:stabbnd}
Consider again the system from Example~\ref{example:stabunbnd} with $A_0=1.5$ and $B_0=1$. Consider that the noise $e(t)$ is drawn uniformly from the set $\{e\,|\, e^2\leq 0.2\}$ and data $(U_-,X)$ is collected for $N=10$. We select four different instrumental variables $r$ based on lagged versions of the input $u$ with $M\in\{2,3,4,5\}$. These are defined as $r_M(t):=\operatorname{col}(u(t),u(t-1),\dots,u(t-M+1)$, i.e., $r_2(t)=\operatorname{col}(u(t),u(t-1))$, $r_3(t)=\operatorname{col}(u(t),u(t-1),u(t-2))$, \emph{et cetera}. We assume prior knowledge on the cross-covariance through the bounds~\eqref{eq:covbndmult} with $c_i^u=-c_i^l=0.1$, $i=1,\dots,M$; these bounds are verified to hold for each of the four choices for $M$. Figure~\ref{fig:lccbmtns} shows the set of feasible systems $\Sigma_{(U_-,X)}^R$ for each choice of $r_M$, denoted $\Sigma_M^R$, illustrating a significant reduction in the size of $\Sigma_M^R$ for increasing $M$. We verify that the data $(U_-,X)$ are informative for quadratic stabilization by Corollary~\ref{cor:quadstab}, since the LMIs~\eqref{eq:infstablin} are feasible for $M=2,\dots,5$, yielding $K=-1.4842$ for $M=5$ such that $\Sigma_K\subseteq \Sigma_5^R$.
\end{example}

\section{Including performance specifications}
In this section, we will consider the problem of finding a feedback gain from the data $(U_-,X)$, such that the closed-loop system with $(A_0,B_0)$ satisfies a given $\mathcal{H}_\infty$ or $\mathcal{H}_2$ performance bound. Consider the performance output $z$, given by
\begin{align*}
z(t)=Cx(t)+Du(t),
\end{align*}
where $C$ and $D$ are user-specified matrices. Recall the set $\Sigma_K$; the set of systems that are stabilized by $K$. The set of systems that achieve $\mathcal{H}_\infty$ performance $\gamma$ with feedback $K$ is defined as
\begin{align*}
\Sigma_{K}^{\mathcal{H}_\infty}(\gamma):=\Sigma_K\cap \{(A,B)\,|\,\|T\|_{\mathcal{H}_\infty}\!\!<\gamma\},
\end{align*}
with $T(q):=C(qI-A-BK)^{-1}+D$.

\begin{definition}
The data $(U_-,X)$ are said to be informative for $\mathcal{H}_\infty$ control with performance $\gamma$ if there exists a feedback gain $K$ so that $\Sigma_{(U_-,X)}^R\subseteq \Sigma_K^{\mathcal{H}_\infty}(\gamma)$.
\end{definition}

\begin{proposition}
Consider a pair $(A,B)$ and $\gamma>0$. The following statements are equivalent:
\begin{itemize}
\item there exists  $K$ so that $(A,B)\in\Sigma_K^{\mathcal{H}_\infty}(\gamma)$,
\item there exist $K$ and $P$ so that $P\succ 0$ and
\begin{align} \label{eq:Hinfeq}
\begin{bmatrix}
I & 0\\ A+BK & I\\ \hline
0 & I\\ C & D
\end{bmatrix}^\top \!\!\left[\begin{array}{cc | cc}
- P & 0 & 0 & 0\\ 0 & P & 0 & 0\\ \hline
0 & 0 & -\gamma^2 I & 0\\ 0 & 0 & 0 & I
\end{array}\right]\!\!\begin{bmatrix}
I & 0\\ A+BK & I\\ \hline
0 & I\\ C & D
\end{bmatrix}\!\prec 0.
\end{align}
\end{itemize}
\end{proposition}
\begin{definition}
The data $(U_-,X)$ are said to be informative for common $\mathcal{H}_\infty$ control with performance $\gamma$ if there exist $K$ and $P$ so that $P\succ 0$ and \eqref{eq:Hinfeq} holds for all $(A,B)\in \Sigma_{(U_-,X)}^R$.
\end{definition}

\begin{proposition}
The data $(U_-,X)$ are informative for common $\mathcal{H}_\infty$ control with performance $\gamma$ if and only if there exist $K$ and $P$ so that $P\succ 0$ and for all $i\in\{1,\dots,L\}$:
\begin{align*}
\begin{bmatrix}
I & 0\\ \sigma_{(U_-,X)}^i\begin{bmatrix}
I\\ K
\end{bmatrix} & I\\ \hline
0 & I\\ C & D%
\end{bmatrix}^{\!\top} \!\!\!\!\!\!\left[\begin{array}{cc | cc}
- P & 0 & 0 & 0\\ 0 & P & 0 & 0\\ \hline
0 & 0 & -\gamma^2 I & 0\\ 0 & 0 & 0 & I
\end{array}\right]\!\!\!\begin{bmatrix}
I & 0\\ \sigma_{(U_-,X)}^i\begin{bmatrix}
I\\ K
\end{bmatrix} & I\\ \hline
0 & I\\ C & D
\end{bmatrix}\!\!\prec 0.
\end{align*}
\end{proposition}

\begin{corollary}
The data $(U_-,X)$ are informative for common $\mathcal{H}_\infty$ control with performance $\gamma$ if and only if there exist $Y$ and $M$ so that for all $i\in\{1,\dots,L\}$:
\begin{align*}
\begin{bmatrix}
Y & 0 & Z^\top (\sigma^i_{(U_-,X)})^\top & YC^\top\\0 & \gamma I & I & D^\top\\
\sigma^i_{(U_-,X)} Z & I & Y & 0\\ CY & D & 0 & \gamma I
\end{bmatrix}\succ 0,
\end{align*}
with $Z:=\operatorname{col} (Y,M)$.
\end{corollary}

The set of systems that achieve $\mathcal{H}_2$ performance $\gamma$ with feedback $K$ is defined as
\begin{align*}
\Sigma_{K}^{\mathcal{H}_2}(\gamma):=\Sigma_K\cap \{(A,B)\,|\,\|T\|_{\mathcal{H}_2}<\gamma\}.
\end{align*}
\begin{definition}
The data $(U_-,X)$ are said to be informative for $\mathcal{H}_2$ control with performance $\gamma$ if there exists a feedback gain $K$ so that $\Sigma_{(U_-,X)}^R\subseteq \Sigma_K^{\mathcal{H}_2}(\gamma)$.
\end{definition}

\begin{proposition}
Consider a pair $(A,B)$ and $\gamma>0$. The following statements are equivalent:
\begin{itemize}
\item there exists  $K$ so that $(A,B)\in\Sigma_K^{\mathcal{H}_2}(\gamma)$,
\item there exist $K$, $P$ and $Z$ so that $\operatorname{trace}Z<\gamma$ and
\begin{align} \label{eq:H2eq}
\!\begin{bmatrix}
P & P(A+BK) & P\\ \star & P & 0\\ \star & \star & \gamma I
\end{bmatrix}\succ 0,\ \begin{bmatrix}
P & 0 & C^\top\\ 0 & I & D^\top \\ C & D & Z
\end{bmatrix}\succ 0.
\end{align}
\end{itemize}
\end{proposition}
\begin{definition}
The data $(U_-,X)$ are said to be informative for common $\mathcal{H}_2$ control with performance $\gamma$ if there exists $K$, $P$ and $Z$ so that $\operatorname{trace} Z<\gamma$ and \eqref{eq:H2eq} holds for all $(A,B)\in \Sigma_{(U_-,X)}^R$.
\end{definition}
\begin{proposition}
The data $(U_-,X)$ are informative for common $\mathcal{H}_2$ control with performance $\gamma$ if and only if there exist $K$, $P$ and $Z$ so that $\operatorname{trace}P<\gamma$ and for all $i\in\{1,\dots,L\}$:
\begin{align*}
\begin{bmatrix}
P & 0 & C^\top\\ 0 & I & D^\top \\ C & D & Z
\end{bmatrix}\succ 0 \quad \text{ and } \quad \begin{bmatrix}
P & P\sigma_{(U_-,X)}^i\begin{bmatrix}
I\\ K
\end{bmatrix} & P\\ \star & P & 0\\ \star & \star & \gamma I
\end{bmatrix}\succ 0.
\end{align*}
\end{proposition}

\begin{corollary}
The data $(U_-,X)$ are informative for common $\mathcal{H}_2$ control with performance $\gamma$ if and only if there exist $Y$, $M$ and $P$ so that $\operatorname{trace}P<\gamma$ and
\begin{align*}
\begin{bmatrix}
Y & 0 & YC^\top\\ 0 & I & D^\top \\ CY & D & P
\end{bmatrix}\succ 0 \quad \text{ and } \quad \begin{bmatrix}
Y & \sigma_{(U_-,X)}^iZ & I\\ \star & Y & 0\\ \star & \star & \gamma I
\end{bmatrix}\succ 0,
\end{align*}
holds for all $i\in\{1,\dots,L\}$ with $Z:=\operatorname{col}(Y,M)$.
\end{corollary}

\section{Concluding remarks}
We have considered the problem of analyzing informativity of data for controller design with prior knowledge on process noise in the form of linear sample cross-covariance bounds. We have established a parametrization of the set of systems that are compatible with data. Using the convexity of this set and the convexity of stability/performance conditions with respect to the system matrices, we have developed necessary and sufficient conditions for informativity for stabilization and $\mathcal{H}_2/\mathcal{H}_\infty$ control. In future work, we will consider the problem of informativity of data for distributed control with cross-covariance bounds.

\small

\bibliographystyle{ifacconf}
\bibliography{../rfrncs21a.bib}

\begin{thebibliography}{16}
\providecommand{\natexlab}[1]{#1}
\providecommand{\url}[1]{\texttt{#1}}
\providecommand{\urlprefix}{URL }
\expandafter\ifx\csname urlstyle\endcsname\relax
  \providecommand{\doi}[1]{doi:\discretionary{}{}{}#1}\else
  \providecommand{\doi}{doi:\discretionary{}{}{}\begingroup
  \urlstyle{rm}\Url}\fi

\bibitem[{Berberich et~al.(2020)Berberich, Koch, Scherer, and
  Allg\"{o}wer}]{berberich2020}
Berberich, J., Koch, A., Scherer, C.W., and Allg\"{o}wer, F. (2020).
\newblock Robust data-driven state-feedback design.
\newblock In \emph{2020 American Control Conf. (ACC)}, 1532--1538.

\bibitem[{Bisoffi et~al.(2021)Bisoffi, {De Persis}, and Tesi}]{bisoffi2021}
Bisoffi, A., {De Persis}, C., and Tesi, P. (2021).
\newblock Trade-offs in learning controllers from noisy data.
\newblock \emph{Systems \& Control Letters}, 154, 104985.

\bibitem[{Campestrini et~al.(2017)Campestrini, Eckhard, Bazanella, and
  Gevers}]{campestrini2017}
Campestrini, L., Eckhard, D., Bazanella, A.S., and Gevers, M. (2017).
\newblock Data-driven model reference control design by prediction error
  identification.
\newblock \emph{J. Franklin Inst.}, 354(6), 2628 -- 2647.

\bibitem[{Campi et~al.(2002)Campi, Lecchini, and Savaresi}]{campi2002}
Campi, M., Lecchini, A., and Savaresi, S. (2002).
\newblock Virtual reference feedback tuning: a direct method for the design of
  feedback controllers.
\newblock \emph{Automatica}, 38(8), 1337 -- 1346.

\bibitem[{Hakvoort and {Van den Hof}(1995)}]{hakvoort95}
Hakvoort, R.G. and {Van den Hof}, P.M.J. (1995).
\newblock Consistent parameter bounding identification for linearly
  parametrized model sets.
\newblock \emph{Automatica}, 31(7), 957--969.

\bibitem[{{Hjalmarsson} et~al.(1998){Hjalmarsson}, {Gevers}, {Gunnarsson}, and
  {Lequin}}]{hjalmarsson98}
{Hjalmarsson}, H., {Gevers}, M., {Gunnarsson}, S., and {Lequin}, O. (1998).
\newblock Iterative feedback tuning: theory and applications.
\newblock \emph{IEEE Control Systems Magazine}, 18(4), 26--41.

\bibitem[{Hou and Wang(2013)}]{hou2013}
Hou, Z.S. and Wang, Z. (2013).
\newblock From model-based control to data-driven control: Survey,
  classification and perspective.
\newblock \emph{Information Sciences}, 235, 3 -- 35.

\bibitem[{Kothare et~al.(1996)Kothare, Balakrishnan, and Morari}]{kothare1996}
Kothare, M.V., Balakrishnan, V., and Morari, M. (1996).
\newblock Robust constrained model predictive control using linear matrix
  inequalities.
\newblock \emph{Automatica}, 32(10), 1361--1379.

\bibitem[{Ljung(1999)}]{ljung1999}
Ljung, L. (1999).
\newblock \emph{System Identification: Theory for the User}.
\newblock Prentice Hall PTR, Upper Saddle River, NJ, USA.

\bibitem[{Scherer and Weiland(2017)}]{schererweilandLMI}
Scherer, C. and Weiland, S. (2017).
\newblock Linear matrix inequalities in control.
\newblock {DISC} lecture notes.

\bibitem[{Steentjes et~al.(2021)Steentjes, Lazar, and {Van den
  Hof}}]{steentjes2021cdc}
Steentjes, T.R.V., Lazar, M., and {Van den Hof}, P.M.J. (2021).
\newblock {$\mathcal{H}_\infty$} performance analysis and distributed
  controller synthesis for interconnected linear systems from noisy input-state
  data.
\newblock In \emph{60\textsuperscript{th} IEEE Conf. Decision and Control
  (CDC)}, 3717--3722. Austin, Texas, USA.

\bibitem[{Steentjes et~al.(2022)Steentjes, Lazar, and Van~den
  Hof}]{steentjes2022qccb}
Steentjes, T.R.V., Lazar, M., and Van~den Hof, P.M.J. (2022).
\newblock On data-driven control: Informativity of noisy input-output data with
  cross-covariance bounds.
\newblock \emph{IEEE Control Systems Letters}, 6, 2192--2197.

\bibitem[{{Van den Hof} and {Schrama}(1995)}]{vandenhofetal1995}
{Van den Hof}, P.M.J. and {Schrama}, R.J.P. (1995).
\newblock Identification and control: Closed-loop issues.
\newblock \emph{Automatica}, 31(12), 1751 -- 1770.

\bibitem[{{van Waarde} et~al.(2020){van Waarde}, {Eising}, {Trentelman}, and
  {Camlibel}}]{waarde2020tac}
{van Waarde}, H.J., {Eising}, J., {Trentelman}, H.L., and {Camlibel}, M.K.
  (2020).
\newblock Data informativity: A new perspective on data-driven analysis and
  control.
\newblock \emph{IEEE Trans. Autom. Control}, 65(11), 4753--4768.

\bibitem[{van Waarde and Camlibel(2021)}]{vanwaarde2021finsler}
van Waarde, H.J. and Camlibel, M.K. (2021).
\newblock A matrix {F}insler's lemma with applications to data-driven control.
\newblock In \emph{60\textsuperscript{th} IEEE Conf. Decision and Control
  (CDC)}, 5770--5775. Austin, Texas, USA.

\bibitem[{van Waarde et~al.(2022)van Waarde, Camlibel, and
  Mesbahi}]{vanwaarde2020noisytac}
van Waarde, H.J., Camlibel, M.K., and Mesbahi, M. (2022).
\newblock From noisy data to feedback controllers: Nonconservative design via a
  matrix {S}-lemma.
\newblock \emph{IEEE Trans. Autom. Control}, 67(1), 162--175.

\end{thebibliography}
\end{document}